\documentclass[12pt]{amsart}
\usepackage{graphicx} 

\usepackage{amsmath,amssymb,amsthm}
\allowdisplaybreaks[1]
\usepackage{thmtools}

\usepackage[backend=bibtex,url=false]{biblatex}
\usepackage{xcolor}

\usepackage{hyperref}
\usepackage{xurl}
\usepackage{cleveref}

\newtheorem{theorem}{Theorem}[section]

\newtheorem{lemma}[theorem]{Lemma}
\newtheorem{corollary}[theorem]{Corollary}
\newtheorem{proposition}[theorem]{Proposition}

\def\N{\mathbb{N}}
\def\Z{\mathbb{Z}}
\def\Q{\mathbb{Q}}

\newcommand{\Bad}{\mathrm{Bad}}
\DeclareMathOperator{\GL}{GL}

\newcommand{\abs}[1]{\left\lvert#1\right\rvert}

\title[The ratio spectrum of Lagrange constants]{The ratio spectrum of Lagrange constants under linear fractional transformations}

\author{Harold Erazo}
\address{SUSTech International Center for Mathematics, Shenzhen, Guangdong, People’s Republic of China}
\email{erazo@sustech.edu.cn}
\thanks{}

\author{Carlos Gustavo Moreira}
\address{Carlos Gustavo Moreira: SUSTech International Center for Mathematics, Shenzhen, Guangdong, People’s Republic of China and IMPA, Rio de Janeiro, Brazil}
\email{gugu@impa.br}
\thanks{The second author was partially supported by CNPq and FAPERJ}

\date{\today}


\begin{document}

\begin{abstract}
In this note we solve a problem posed by Lagarias and Shallit concerning Lagrange constants under linear fractional transformations $Mx=\frac{ax+b}{cx+d}$. For an integer matrix $M$ with nonzero determinant and relatively prime entries, define the ratio spectrum
\begin{equation*}
    \mathcal{V}(M)=\left\{\frac{k(Mx)}{k(x)}:x\in\Bad\right\},
\end{equation*}
where $k(x)$ denotes the Lagrange constant of the irrational number $x$ and $\Bad$ is the set of badly approximable numbers. Lagarias and Shallit proved that
\begin{equation*}
    \mathcal{V}(M)\subseteq\left[\frac{1}{|\det M|},|\det M|\right],
\end{equation*}
and asked for the determination of $\mathcal{V}(M)$. We prove that
\begin{equation*}
    \mathcal{V}(M)=\left[\frac{1}{|\det M|},|\det M|\right].
\end{equation*}
\end{abstract}

\maketitle

\section{Introduction}

Given an irrational $x$, we denote by $x=[a_0;a_1,\dots]$ its simple continued fraction and let $K_\infty(x):=\limsup_{i\geq 1}a_i$. It is well known that $x$ is \emph{badly approximable} if and only if $K_\infty(x)<\infty$. We denote by $\Bad$ the set of badly approximable numbers.

We use $M_2(\Z)$ to denote the set of $2\times2$ integer matrices and $\GL_2(\Z)$ to denote the set of $2\times2$ matrices with determinant $\pm 1$. Given an integer matrix
\begin{equation*}
    M = \begin{pmatrix}
        a & b \\
        c & d
    \end{pmatrix}\in M_2(\Z), \qquad \det M = ad-bc\neq 0,
\end{equation*}
and a real number $x$, we let whenever possible
\begin{equation*} 
    Mx=\frac{ax+b}{cx+d}. 
\end{equation*}

Since this action is invariant by scalar multiplication to $M$, throughout the paper, unless otherwise stated, all matrices have nonzero determinant and $\gcd(a,b,c,d)=1$. 

A classical theorem states that two irrational numbers have eventually equal continued fraction expansions if and only if they differ by a transformation in $\GL_2(\Z)$. A more general fact, that follows for example from \cite{Raney}, is that for every matrix $M\in M_2(\Z)$ with $\det M\neq 0$, the irrational number $Mx$ has bounded partial quotients if and only if $x$ has bounded partial quotients.

This result was quantitatively improved by Lagarias and Shallit in \cite{LagariasShallit, LagariasShallitCorrection} and more recently in \cite{Stambul, Turkish}.  
\begin{theorem}[Lagarias-Shallit]\label{thm:LagariasShallitKinfty}
Let $M\in M_2(\Z)$ with $\det M\neq 0$. Given any $x\in\Bad$ we have
\begin{equation*}
    \frac{1}{|\det M|}K_{\infty}(x)-2\leq K_{\infty}(Mx)\leq|\det M|(K_\infty(x)+2).
\end{equation*}
\end{theorem}

Earlier, Davenport \cite{Davenport} showed that for any irrational $x$ and prime $p$, one has $K_\infty(x+\frac{a}{p})\geq p-1$ for some $0\leq a<p$. Mendès France \cite{MendesImprovement} refined this result by showing that for some $0\leq a<p$, the continued fraction $x+\frac{a}{p}=[0;b_1,b_2,\dots]$ satisfies $b_n\geq p-1$ for a positive proportion of $n$.

Given an irrational $x=[a_0;a_1,a_2,\dots]$, we define its \emph{Lagrange constant} as 
\begin{align*}
    k(x) &:= \limsup_{p,q\to\infty} (|q(qx - p)|)^{-1} \\
    &= \sup\left\{k>0\,:\,\left|x - \frac{p}{q}\right| < \frac{1}{kq^{2}} \text{ has infinitely many solutions }\frac{p}{q}\in\Q\right\}.
\end{align*}

Since good approximations are precisely the convergents because of Legendre's theorem, this quantity can be computed in terms of the continued fraction of $x$ as 
\begin{equation}\label{eq:k(x)}
    k(x)=\limsup_{n\to\infty}(q_n|q_nx-p_n|)^{-1}=\limsup_{n\to\infty}([a_n;a_{n+1},\dots]+[0;a_{n-1},\dots,a_1]),
\end{equation}
where $p_n/q_n=[a_0;a_1,a_2,\dots,a_n]$ are the convergents of $x$.

Based on a result by Cusick and Mendès \cite{CusickMendes}, Lagarias and Shalit also proved quantitative estimates for the Lagrange constant \cite[Theorem 3.2]{LagariasShallit} of a badly approximable number after applying a linear fractional transformation.

\begin{theorem}[Lagarias-Shallit]\label{thm:Lagarias-Shallit_Lagrange_constant}
Let $M\in M_2(\Z)$ with $\det M\neq 0$. Given any $x\in\Bad$ we have
\begin{equation*}
    \frac{1}{|\det M|}\cdot k(x)\leq k(Mx)\leq |\det M|\cdot k(x).
\end{equation*}    
\end{theorem}

We shall give a short proof of \Cref{thm:Lagarias-Shallit_Lagrange_constant} based on Smith normal form. We remark that the same argument extends to other principal ideal domains, such as the rings of integers $\mathcal{O}_{\mathbb{K}}$ of the imaginary quadratic fields $\mathbb{K}=\Q(\sqrt{-D})$, where $D\in\{1, 2, 3, 7, 11, 19, 43, 67, 163\}$, for matrices $M\in M_2(\mathcal{O}_{\mathbb{K}})$ with $\det M\neq 0$ and the corresponding approximation constants defined with respect to $p,q\in \mathcal{O}_{\mathbb{K}}$.

Lagarias and Shallit also observed that both bounds are sharp. More precisely, for every $M\in M_2(\Z)$ with $\det M\neq 0$, there exist irrational numbers $x$ and $y$ such that $k(Mx)/k(x)=|\det M|$ and $k(My)/k(y)=1/|\det M|$. 

This led them to ask \cite[Page 278]{LagariasShallit} for a determination of the set
\begin{equation}
    \mathcal{V}(M):=\left\{\frac{k(Mx)}{k(x)}:x\in\Bad\right\}.
\end{equation}

Using Smith normal form, it is not difficult to reduce this problem to the case of multiplication by a positive integer. More precisely, it suffices to consider matrices of the form $M=\begin{pmatrix}
    m & 0 \\
    0 & 1
\end{pmatrix}$ in which case $Mx=mx$. Using ideas inspired by Hall's construction \cite{Hall}, we prove that the ratio spectrum coincides with the entire interval allowed by \Cref{thm:Lagarias-Shallit_Lagrange_constant}.

\begin{theorem}\label{thm:main}
Given any $M=\begin{pmatrix} a & b \\ c & d \end{pmatrix}\in M_2(\Z)$, with $\det M\neq 0$ and $\gcd(a,b,c,d)=1$, we have that 
\begin{equation*}
    \mathcal{V}(M)=\left[\frac{1}{|\det M|},|\det M|\right].
\end{equation*}
\end{theorem}

Moreover, it possible to show that for any $|\det M|^{-1}\leq \gamma\leq |\det M|$, the set of $x\in\Bad$ such that $k(Mx)/k(x)=\gamma$ has Hausdorff dimension equal to 1.

From an algorithmic perspective, several methods have been developed to compute the continued fraction expansion of a linear fractional transformation $Mx$ from that of $x$. In 1894, Hurwitz \cite{Hurwitz1963} gave an algorithm to obtain the continued fraction of $2x$ from that of $x$. Hall \cite{Hall} later introduced a general procedure for computing the continued fraction expansion of $Mx$. In \cite{Raney}, Raney introduced a much simpler method by introducing transducers. Other algorithms have also been proposed, including Gosper's algorithm \cite[Page 649]{Knuth}. 

More recently, several authors have investigated the relation between the Lagrange constants of $x$ and $Mx$. Results in this direction can be found in \cite{Einsiedler, Compositio, BFS}.

\section{Proof of Theorem \ref{thm:main}}

We need some basic facts about continued fractions. We have the following exact formula 
\begin{equation}\label{eq:alpha_n_and_beta_n}
    \abs{\alpha-\frac{p_n}{q_n}} =\frac1{(\alpha_{n+1}+\beta_{n+1}) q_n^2}, \qquad \text{for all } n\in {\mathbb N},
\end{equation}
where $p_n/q_n=[a_0;a_1,a_2,\dots,a_n]$, $\alpha_{n}=[a_n;a_{n+1},a_{n+2},\ldots]$ and $\beta_n=[0;a_{n-1},a_{n-2},\ldots,a_1]$. In particular

\begin{lemma}\label{lem:approx_convg}
    Let $x=[a_0,a_1,a_2,\dots]$ and $\frac{p_n}{q_n}=[a_0,a_1,\dots, a_n]$. Then
    \begin{equation}\label{eq:coefficients_inequality}
        \frac{1}{(a_{n+1}+2)q_n^2}<\abs{x-\frac{p_{n}}{q_{n}}}<\frac{1}{a_{n+1}q_n^2}.
    \end{equation}
\end{lemma}

\begin{lemma}[Legendre]\label{lem:legendre}
    If 
    \begin{equation*}
        \abs{x - \frac{p}{q}}<\frac{1}{2q^2}, \qquad \frac{p}{q}\in\Q,
    \end{equation*}
    then $\frac{p}{q}$ is a convergent of $x$ (i.e., $\frac{p}{q}=\frac{p_n}{q_n}$ for some $n\in\mathbb N$).
\end{lemma}

Another classical ingredient is Hall's theorem on sums of Gauss--Cantor sets \cite{Hall}. For $n\geq 1$, let
\begin{equation*}
    \mathcal{C}(n) = \{[0;a_1,a_2,\dots]:1\leq a_i\leq n \text{ for all } i\geq 1\}.
\end{equation*}

For every $n\geq 4$, we have
\begin{equation}\label{eq:Hall}
    \mathcal{C}(n)+\mathcal{C}(n)=[2\min C(n),2\max C(n)].
\end{equation}
It is easy to check that $2(\max \mathcal{C}(n) - \min \mathcal{C}(n)) > 1$ for $n\geq 4$.

The first step to prove \Cref{thm:main} is to reduce the problem to the multiplication by a fixed integer using its Smith normal form \cite[Chapter VII, Proposition 2.11]{Hungerford}. Given any matrix with integer coefficients $M\in M_2(\Z)$, we can find matrices $U,V\in \GL_2(\Z)$ and a diagonal matrix $D^\prime\in M_2(\Z)$ with positive entries
\begin{equation*}
    D^\prime = \begin{pmatrix}
        m & 0 \\
        0 & d
    \end{pmatrix},
\end{equation*}
such that $d$ divides $m$ and such that $M=UD^\prime V$. Since the entries of $M$ are relatively prime, we have that in fact $d=1$. In particular, writing $y=Vx$, we will have that $k(x)=k(y)$ so 
\begin{equation}\label{eq:k(my)k(y)}
    \frac{k(Mx)}{k(x)}=\frac{k(D^\prime Vx)}{k(x)}=\frac{k(D^\prime y)}{k(y)}=\frac{k(my)}{k(y)}.
\end{equation}

\begin{proof}[Proof of \Cref{thm:Lagarias-Shallit_Lagrange_constant}]
    Let $x\in\Bad$ and $M\in M_2(\Z)$ with $\det M\neq 0$. By \eqref{eq:k(my)k(y)} we have that there is an irrational $y$ and a positive integer $m$ such that $k(Mx)/k(x)=k(my)/k(y)$. 

    Since $k(y)=\limsup_{q\to\infty}(q\cdot\lVert qy\rVert)^{-1}$ where $\lVert\cdot\rVert$ denotes the distance to the nearest integer, we have
    \begin{align*}
        \frac{1}{m}\cdot\limsup_{q\to\infty}(q\cdot\lVert qy\rVert)^{-1}&\leq \limsup_{q\to\infty}(q\cdot\lVert q(my)\rVert)^{-1} \\
        &=m\cdot\limsup_{q\to\infty}((mq)\cdot\lVert (mq)y\rVert)^{-1} \\
        &\leq m\cdot\limsup_{q\to\infty}(q\cdot\lVert qy\rVert)^{-1}
    \end{align*}
    which is precisely $\frac{1}{m}\cdot k(y)\leq k(my)\leq m\cdot k(y)$.

    \end{proof}

    Given a finite word of positive integers $w=(a_1,\dots,a_n)\in(\N_{>0})^{n}$, we can define the matrix
    \begin{equation}\label{eq:matrix_M_w}
        M_w=\begin{pmatrix}
            a_1 & 1 \\
            1 & 0
        \end{pmatrix}
        \dotsb
        \begin{pmatrix}
            a_n & 1 \\
            1 & 0
        \end{pmatrix}, \quad\text{so}\quad 
        \begin{pmatrix}
            a_0 & 1 \\
            1 & 0
        \end{pmatrix}
        M_w=
        \begin{pmatrix}
        p_n & p_{n-1} \\
        q_n & q_{n-1}
    \end{pmatrix}.
    \end{equation}

    It is well known that the matrices
    \begin{equation*}
        S = 
        \begin{pmatrix} 
            0 & 1 \\ 
            1 & 0 
        \end{pmatrix} 
        \qquad\text{and}\qquad 
        T = 
        \begin{pmatrix} 
        1 & 1 \\ 
        1 & 0 
        \end{pmatrix}
    \end{equation*}
    generate $GL_2(\Z)$ as a group. We need the following simple corollary:
    \begin{lemma}
        Let $m\geq 2$ be a positive integer. The matrices
        \begin{equation*}
        \begin{pmatrix} 
            m & 1 \\ 
            1 & 0 
        \end{pmatrix} 
        \qquad\text{and}\qquad 
        \begin{pmatrix} 
        1 & 1 \\ 
        1 & 0 
        \end{pmatrix},
    \end{equation*}
        generate $\{A\in\GL_2(\Z/m\Z):\det A=\pm 1 \pmod{m}\}$ as a semigroup.
    \end{lemma}
    \begin{proof}
        Notice that $\GL_2(\Z/m\Z)$ is a finite group. Hence, every semigroup contained in it is a group. Modulo $m$, the two matrices above reduce to $S$ and $T$, respectively. Since $S$ and $T$ generate $\GL_2(\Z)$ as a group, their reductions generate the image of $\GL_2(\Z)$ in $\GL_2(\Z/m\Z)$.
    \end{proof}
    
    \begin{corollary}\label{cor}
        Let $m\geq 2$ be a positive integer. There is a constant $C_m$ depending only on $m$ such that the following holds: given finite words $w_1=(a_1,\dots,a_{n_1})\in(\N_{>0})^{n_1}$ and $w_2=(b_1,\dots,b_{n_2})\in(\N_{>0})^{n_2}$ and a matrix $M^\prime\in\GL_2(\Z)$, then there is a positive integer $1\leq\ell\leq C_m$ and a finite word $(c_1,\dots,c_\ell)\in\{1,\dots,m\}^{\ell}$, such that
        \begin{equation*}
            M_{w_1}\begin{pmatrix}
            c_1 & 1 \\
            1 & 0
        \end{pmatrix}
        \dotsb
        \begin{pmatrix}
            c_\ell & 1 \\
            1 & 0
        \end{pmatrix} M_{w_2} \equiv
        M^\prime \pmod{m}
        \end{equation*}
        where the congruences are entrywise. 
    \end{corollary}

    We will obtain \Cref{thm:main} as a corollary of a stronger proposition.

    \begin{proposition}\label{prop:main}
        Let $m\geq 2$ be a positive integer. Let $C, D \geq \max\{m+1+\sqrt{2},5+\sqrt{2}\}$ be constants such that $1/m^2 \le C/D \le 1$. Then there is an irrational number $y$ such that $k(y)=D$ and $k(my)=mC$.
    \end{proposition}

    \begin{proof}
    Choose the positive integers $A=\lfloor C-(\sqrt{2}-1)\rfloor$ and $B=\lfloor D-(\sqrt{2}-1)\rfloor$ where $\lfloor\cdot\rfloor$ denotes the integer part. Then $A, B \geq \max\{m+2,6\}$. By \eqref{eq:Hall}, there exist $\xi_1,\xi_2,\eta_1,\eta_2\in \mathcal{C}(4)$ such that
    \begin{equation*}
        A+\xi_1+\xi_2=C, \qquad B+\eta_1+\eta_2=D.
    \end{equation*}
    Write
    \begin{equation*}
        \xi_1=[0;a_1,a_2,\ldots],\qquad
        \xi_2=[0;a_{-1},a_{-2},\ldots],
    \end{equation*}
    \begin{equation*}
        \eta_1=[0;b_1,b_2,\ldots],\qquad
        \eta_2=[0;b_{-1},b_{-2},\ldots],
    \end{equation*}
    where $1\leq a_i\leq 4$ and $1\leq b_i\leq 4$ for all $i\neq 1$.

    Now we construct an irrational number $y=[0;d_1,d_2,\dots]$ satisfying $k(y)=D$ and $k(my)=mC$. To begin, let $(d_1,d_2,d_3)=(a_{-1},A,a_1)$. By \Cref{cor}, there exists a word $(c_1^{(1)},\dots,c_{\ell_1}^{(1)})\in\{1,\dots,m\}^{\ell_1}$ with $1\leq\ell_1\leq C_m$ such that 
    \begin{equation*}
        \begin{pmatrix}
            0 & 1 \\
            1 & 0
        \end{pmatrix}
        M_{(a_{-1},A,a_1)}\begin{pmatrix}
        c_1^{(1)} & 1 \\
        1 & 0
    \end{pmatrix}
    \dotsb
    \begin{pmatrix}
        c_{\ell_1}^{(1)} & 1 \\
        1 & 0
    \end{pmatrix}
    \begin{pmatrix}
        b_{-1} & 1 \\
        1 & 0
    \end{pmatrix}
    \equiv
    \begin{pmatrix}
        1 & 1 \\
        1 & 0
    \end{pmatrix} \pmod{m}
    \end{equation*}
    We then define $(d_{4},\dots,d_{\ell_1+6})=(c_1^{(1)},\dots,c_{\ell_1}^{(1)},b_{-1},B,b_1)$. Applying \Cref{cor} once more, we choose $(c_1^{(2)},\dots,c_{\ell_2}^{(2)})\in\{1,\dots,m\}^{\ell_2}$ with $1\leq\ell_2\leq C_m$ such that 
    \begin{equation*}
    \begin{pmatrix}
        0 & 1 \\
        1 & 0
    \end{pmatrix}
        M_{(d_1,\dots,d_{\ell_1+6})}\begin{pmatrix}
        c_1^{(2)} & 1 \\
        1 & 0
    \end{pmatrix}
    \dotsb
    \begin{pmatrix}
        c_{\ell_2}^{(2)} & 1 \\
        1 & 0
    \end{pmatrix}
    M_{(a_{-2},a_{-1})}
    \equiv
    \begin{pmatrix}
        1 & 0 \\
        0 & 1
    \end{pmatrix} \pmod{m}
    \end{equation*}

    Proceeding inductively and alternating between the target matrices $T$ and $I$, we obtain an irrational number
        \begin{equation*}
        y = [0;d_1,d_2,\dots] = [0;a^{(1)},c^{(1)},b^{(1)},c^{(2)},a^{(2)},c^{(3)},b^{(2)},\dots]
    \end{equation*}
    where we have denoted $a^{(k)}=(a_{-k},\dots,a_{-1},A,a_1,\dots,a_k)$, \\ $b^{(k)}=(b_{-k},\dots,b_{-1},B,b_1,\dots,b_k)$ and $c^{(k)}=(c_1^{(k)},\dots,c_{\ell_k}^{(k)})$, which by construction satisfies
    \begin{itemize}
        \item $d_n\leq \max\{m,4\}$ whenever $d_n\not\in\{A,B\}$,
        \item $q_{n}\equiv 0\pmod{m}$ whenever $d_{n+1} = A$,
        \item $q_{n}\equiv 1\pmod{m}$ whenever $d_{n+1} = B$.
    \end{itemize}

    Since $A$ and $B$ are greater than $m+2$ and we have $B+\eta_1+\eta_2 = D\geq C = A+\xi_1+\xi_2$, by \eqref{eq:k(x)} we have that the Lagrange constant of $y$ is equal 
    \begin{align}
        k(y)&=\limsup_{n\to\infty}([d_n;d_{n+1},\dots]+[0;d_{n-1},\dots,d_1]) \nonumber \\
        &=\limsup_{n\to\infty}([B;b_{1},\dots,b_n,\dots]+[0;b_{-1},\dots,b_{-n}]) = B+\eta_1+\eta_2 = D \label{eq:k(y)}
    \end{align}

    On the other hand, if $p/q\in\Q$ with $\gcd(p,q)=1$ is a very good approximation of $my$, namely, if
    \begin{equation*}
        \abs{my-\frac{p}{q}}<\frac{1}{m\cdot\max\{m+2,6\}\cdot q^2},
    \end{equation*}
    then 
    \begin{equation*}
        \abs{y-\frac{p}{mq}}<\frac{1}{\max\{m+2,6\}\cdot(mq)^2},
    \end{equation*}
    whence by \Cref{lem:legendre} we have that $p/mq$ is a convergent of $y$ (this trick was already noticed by Perron \cite[Chapter IV]{Perron2013}). Moreover, by \eqref{eq:coefficients_inequality}, such convergent $p_{n}/q_{n}=p/mq$ must satisfy that $d_{n+1}=A$ or $d_{n+1}=B$. In particular, we have that $(p,q)=(p_n,q_n/m)$ when $d_{n+1}=A$ and $(p,q)=(mp_n,q_n)$ when $d_{n+1}=B$. Therefore, using the notation of \eqref{eq:alpha_n_and_beta_n}, we have
    \begin{align*}
        k(my) &= \limsup_{p,q\to\infty} \frac{1}{|q(q(my) - p)|} \\
        &= \max\left\{\limsup_{\substack{n\to\infty \\ d_{n+1}=A}} \frac{1}{|(q_n/m)((q_n/m)(my) - p_n)|},\limsup_{\substack{n\to\infty \\ d_{n+1}=B}}\frac{1}{|q_n(q_n(my) - mp_n)|}\right\} \\
        &=\max\left\{\limsup_{\substack{n\to\infty \\ d_{n+1}=A}}m(\alpha_{n+1}+\beta_{n+1}),\limsup_{\substack{n\to\infty \\ d_{n+1}=B}}\frac{\alpha_{n+1}+\beta_{n+1}}{m}\right\} \\
        &=\max\left\{m(A+\xi_1+\xi_2),\frac{B+\eta_1+\eta_2}{m}\right\} \\
        &=\max\left\{mC,\frac{D}{m}\right\} =mC
    \end{align*}
    since $1/m^2\leq C/D$ by hypothesis, which finishes the proof of \Cref{prop:main}.
    \end{proof}

    As a final remark, it is not difficult to construct a set of irrational numbers of Hausdorff dimension $1$ satisfying $k(mz)/k(z)=mC/D$. Indeed, let $(n_j)_{j\geq 1}$ be a sufficiently rapidly growing sequence, for instance $n_j=2^j$. Choose integers $A,B$ and $\xi_1,\xi_2,\eta_1,\eta_2\in \mathcal{C}(4)$ as above. Consider all irrational numbers of the form
    \begin{equation*}
        z=[0;e_1,\dots,e_{n_1},c^{(0)},a^{(1)},c^{(1)},b^{(1)},e_{n_1+1},\dots,e_{n_1+n_2},c^{(2)},a^{(2)},c^{(3)},b^{(2)},\dots]
    \end{equation*}
    where $1\leq e_i \leq A-2$ for all $i\geq 1$. Here the finite words $c^{(k)}$ also depend on the coefficients $e_i$, but it is irrelevant because their length is always uniformly bounded by $C_m$. Adapting the argument from the proof of \cite[Theorem 2]{geometricproperties}, one obtains that the Hausdorff dimension of this set is comparable to $\dim_H(\mathcal{C}(A-2))$, whose value tends to $1$ as $A\to\infty$.

    \section*{Acknowledgements}
    This work was partially supported by the Simons Foundation grant (award no. SFI-MPS-T-Institutes-00010825) and from State Treasury funds as part of a task commissioned by the Minister of Science and Higher Education under the project “Organization of the Simons Semesters at the Banach Center - New Energies in 2026-2028” (agreement no. MNiSW/2025/DAP/491). We also thank the referee for their thorough and careful reading of the paper. 

\vspace{0.5cm} 
\begin{center}
    \includegraphics[width=0.8\textwidth]{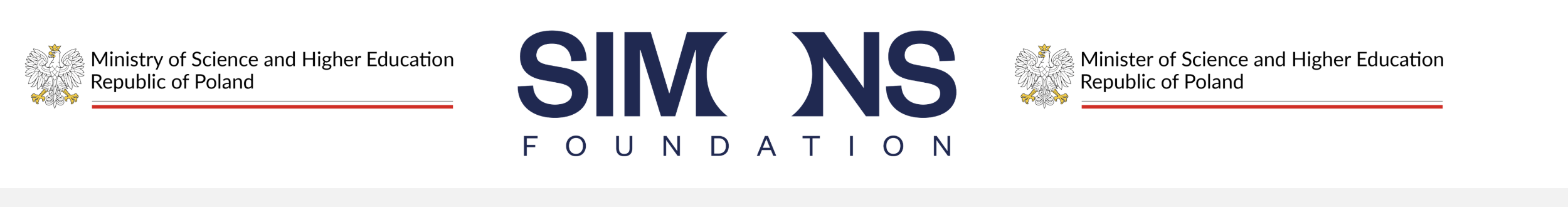}
\end{center}

\printbibliography

\end{document}